\documentclass[12pt]{article}
\usepackage{amsfonts,amsmath,amssymb}
\usepackage{enumerate}
\usepackage{theorem}

\numberwithin{equation}{section}

\newtheorem{theorem}{Theorem}[section]
\newtheorem{lemma}[theorem]{Lemma}
\newtheorem{proposition}[theorem]{Proposition}

\newtheorem{definition}{Definition}[section]
\newtheorem{corollary}[theorem]{Corollary}

\newcommand{\cl}[1]{\mathcal{#1}} 
\newcommand{\bb}[1]{\mathbb{#1}}

\newcommand{\sca}[1]{\left\langle#1\right\rangle} 
\newcommand{\nor}[1]{\left\Vert #1\right\Vert}

\begin{document}

\title{A MORITA TYPE EQUIVALENCE FOR DUAL OPERATOR ALGEBRAS}
\author{G.K. ELEFTHERAKIS}

\date{}

\maketitle

\begin{abstract}
We generalize the main theorem of Rieffel for Morita equivalence of
$W^*$-algebras to the case of unital dual operator algebras: two unital
dual operator algebras $\cl{A}, \cl{B}$ have completely isometric normal
representations $\alpha, \beta $ such that
$\alpha(\cl{A})=[\cl{M}^*\beta(\cl{B})\cl{M}]^{-w^*}$ and
$\beta(\cl{B})=[\cl{M}\alpha(\cl{A})\cl{M}^*]^{-w^*}$ for a ternary ring of operators $\cl{M}$
(i.e. a linear space $\cl{M}$ such that $\cl{M}\cl{M}^*\cl{M}\subset \cl{M})$ if and only if there exists an equivalence  functor
$\cl{F}:\; _{\cl{A}}\mathfrak{M}\rightarrow \; _{\cl{B}}\mathfrak{M}$
which ``extends" to a $*-$functor implementing an equivalence between the categories
$_{\cl{A}}\mathfrak{DM}$ and $_{\cl{B}}\mathfrak{DM}.$ By $_{\cl{A}}\mathfrak{M}$
we denote the category of normal representations of $\cl{A}$ and by $_{\cl{A}}\mathfrak{DM}$
the category with the same objects as $_{\cl{A}}\mathfrak{M}$ and $\Delta (\cl{A})$-module
maps as morphisms ($\Delta (\cl{A})=\cl{A}\cap \cl{A}^*$).
We prove that this functor is equivalent to a functor ``generated" by a
$\cl{B}, \cl{A}$ bimodule, that it is normal and completely isometric.
\end{abstract}



\section{Introduction}

In the beginning of the 70's, M. Rieffel \cite{rif} (see also \cite{rif2}) introduced
to Operator Theory the notion of Morita equivalence. Rieffel's work was concerned with the
equivalence of representations of $C^*$ and $W^*$ algebras.
With the development of the theory of operator spaces,
it was natural to seek extensions of this theory to the class of (abstract) operator algebras.

The papers \cite{bmp} and \cite{bmor} deal with Morita equivalence of not necessarily selfadjoint (norm closed) operator algebras. To this day however, as far as we know, there is no complete theory of
Morita equivalence for dual operator algebras.
A natural  requirement for such a theory would be to respect the additional topological
structure that dual operator algebras possess as dual operator spaces.
 A step in this direction is taken in \cite{bmag}, where Rieffel's theory of Hilbert
 modules is extended to (dual) modules over dual (nonselfadjoint) operator algebras. In this paper
 we are able to generalize Rieffel's theory in a different direction. We study a new notion of
equivalence for representations of dual operator algebras on Hilbert spaces. This equivalence
coincides in the $W^*$-algebra case with the one studied by M. Rieffel; in the non selfadjoint
case there are differences in that two distinct categories have to be simultaneously equivalent.
 We will say that two unital dual operator algebras are {\em $\Delta$-equivalent}
when there is an equivalence functor between their normal representations which not only
preserves intertwiners of representations of the algebras, but also preserves intertwiners of
restrictions to  the diagonals (see Definition \ref{1.4.d}). In \cite{ele} a new notion of equivalence between concrete $w^*$ closed operator algebras 
was developed:
\begin{definition}\label{1.1.d}\cite{ele}  Let $\cl{A}, \cl{B}$ be $w^*$ closed algebras 
acting on Hilbert spaces $H_1$ and $H_2$ respectively. If there is a TRO 
$\cl{M}\subset B(H_1,H_2)$ (i.e. a subspace of $B(H_1,H_2)$ satisfying 
$\cl{M}\cl{M}^*\cl{M}\subset \cl{M})$ ) such that $\cl{A}=[\cl{M}^*\cl{B}\cl{M}]^{-w^*}\;\; \text{and}\;\;
\cl{B}=[\cl{M}\cl{A}\cl{M}^*]^{-w^*}$ we write $\cl{A} 
\stackrel{\cl{M}}{\sim}\cl{B}.$ The algebras $\cl{A}, \cl{B}$
 are called \textbf{TRO equivalent} if there is a TRO 
$\cl{M}$ such that $\cl{A} 
\stackrel{\cl{M}}{\sim}\cl{B}.$ 
\end{definition}
Our first main theorem (Theorem \ref{1.6}) which generalizes the main result in \cite{rif}
is that two (abstract) unital
dual operator algebras  $\cl{A}, \cl{B}$ are
$\Delta $-equivalent if  and only if they have  completely isometric normal representations
$\alpha, \beta $ such that the algebras $\alpha (\cl{A}), \beta (\cl{B})$ are TRO equivalent.
The second main theorem (Theorem \ref{6.3}) states that every $\Delta $-equivalent functor is (unitarily) equivalent 
to a functor ``generated" by an algebra bimodule. The bimodule is generated by ``saturating" 
the TRO which implements the equivalence.

 We present some symbols used below. If $\cl{A}$ is an operator algebra 
we denote its diagonal $\cl{A}\cap \cl{A}^*$ by $\Delta (\cl{A}).$ 
The symbol $[\cl{S}]$ denotes the linear span of
$\cl{S}.$ The commutant of a set $\cl{L}$ of bounded operators on a 
Hilbert space $H$ is denoted $\cl{L}^\prime.$ If $\cl{U}$ 
is a linear space and $n,m\in \mathbb{N}$ we denote by $M_{n,m}(\cl{U})$ the space of
 $n\times m$ matrices with entries from $\cl{U}$ and by $M_{n}(\cl{U})$ the space 
$M_{n,n}(\cl{U}).$ If $\cl{U}, \cl{V}$ are linear spaces, $\alpha $ is  a linear map 
from $\cl{U}$ to $\cl{V}$ and $n,m\in \mathbb{N}$ we denote the linear map\vspace*{-1ex} 
$$M_{n,m}(\cl{U})\rightarrow  M_{n,m}(\cl{V}): (A_{ij})_{i,j}\rightarrow 
(\alpha (A_{ij}))_{i,j}$$ 
again by $\alpha .$ If $\cl{U}$ is a subspace of 
$B(H,K)$ for $H,K$ Hilbert spaces we equip $M_{n,m}(\cl{U}) $, $n,m\in \mathbb{N}$ 
with the norm inherited from the embedding $M_{n,m}(\cl{U})\subset B(H^n,K^m).$ If 
$(\cl{X}, \|\cdot\|)$ is a normed space we denote by $Ball(\cl{X})$ the unit ball of
 $\cl{X}:$ $\{X\in \cl{X}: \|X\| \leq 1 \}.$ If $x_1,...,x_n$ are in a vector space 
$\cl{V},$ we write $(x_1,...,x_n)^t$ for the column vector in $M_{n,1}(\cl{V}). $

 We present some definitions and concepts used in this work. 
A $C^*$ algebra which is a dual Banach space is called a \textbf{$W^*$ algebra}. 
A \textbf{dual operator algebra} is an operator algebra which is the dual of an 
operator space. Every $W^*$ algebra is a dual operator algebra. For 
every dual operator algebra $\cl{A}$ there exists a Hilbert space $H_0$ and an 
algebraic homomorphism $\alpha _0: \cl{A}\rightarrow B(H_0)$ which is a complete 
isometry and a $w^*$-continuous map \cite{bm}.
 
\begin{lemma}\label{1.1}\cite[8.5.32]{bm}  Let $\cl{C}, \cl{E}$ be von Neumann algebras 
acting on Hilbert spaces $H_1$ and $H_2$ respectively, $\theta : \cl{C}\rightarrow \cl{E}$ 
 be a $*$-isomorphism and
$$\cl{M}=\{T\in B(H_1,H_2): TA=\theta (A)T \; \text{for\;\;all}\;\; A\in \cl{C}\}.$$ 
Then the space $\cl{M}$ is an essential TRO, i.e. the algebras $ [\cl{M}^*\cl{M}]^{-w^*}, 
[\cl{M}\cl{M}^*]^{-w^*}$ contain the identity operators.   
\end{lemma}

We now define the category $_{\cl{A}}\mathfrak{M}$ for a unital dual operator
algebra $\cl{A}$ \cite{bm}. The objects of $_{\cl{A}}\mathfrak{M}$ are pairs 
$(H,\alpha)$ where $H$ is a Hilbert space and 
$\alpha: \cl{A}\rightarrow B(H)$ is a \textbf{normal representation} of $\cl{A}$, i.e.
 a unital completely contractive $w^*$-continuous homomorphism.
If $(H_i,\alpha_i ), i=1,2$ are objects of the category $_{\cl{A}}\mathfrak{M}$ 
the space of homomorphisms $ \mathrm{ Hom}_{\cl{A}}(H_1,H_2)$ is the following: 
\vspace*{-1ex}
$$\mathrm{ Hom}_{\cl{A}}(H_1,H_2)=\{T\in B(H_1,H_2): T\alpha_1(A)= \alpha_2(A)T \text{\;for\;all\;}A\in \cl{A}\}.$$
Observe that the map $\alpha _i|_{\Delta (\cl{A})}$ is a $*$-homomorphism since
$\alpha _i$ is a contraction \cite{bm}. We also define the category  $_{\cl{A}}\mathfrak{DM}$ which has the same objects as
 $_{\cl{A}}\mathfrak{M}$ but for every pair of objects $(H_i,\alpha_i ), i=1,2$ the space 
of homomorphisms $\mathrm{ Hom}_{\cl{A}}^{\mathfrak{D}}(H_1,H_2)$ is given by: 
\vspace*{-1ex}
$$\mathrm{ Hom}_{\cl{A}}^{\mathfrak{D}}(H_1,H_2)=\{T\in B(H_1,H_2): T\alpha_1(A)= 
\alpha_2(A)T \text{\;for\;all\;}A\in \Delta (\cl{A})\}.$$
 If $\cl{A}$ is a $W^*$-algebra the categories 
$_{\cl{A}}\mathfrak{M}$ and $_{\cl{A}}\mathfrak{DM}$ are the same. Also observe 
that $\mathrm{ Hom}_{\cl{A}}(H_1,H_2)\subset \mathrm{ Hom}_{\cl{A}}^{\mathfrak{D}}(H_1,H_2).$

\begin{definition}\label{1.2.d}Let $\cl{A}, \cl{B}$ be unital dual operator 
algebras and $\cl{F}:\;_{\cl{A}}\mathfrak{M} \rightarrow\; _{\cl{B}}\mathfrak{M}$ be a 
functor.We say that the functor  $\cl{F}$ has a $\Delta$-extension if there is a 
functor $\cl{G}: \;_{\cl{A}}\mathfrak{DM} \rightarrow \;_{\cl{B}}\mathfrak{DM}$ 
such that the following diagram is commutative:
\begin{align*}\begin{array}{clr} _{\cl{A}}\mathfrak{M}  & \hookrightarrow  
& _{\cl{A}}\mathfrak{DM} \\  \cl{F}\downarrow &  & \cl{G}\downarrow\;\;\;\;\;\;  
\\  _{\cl{B}}\mathfrak{M}  & \hookrightarrow  & _{\cl{B}}\mathfrak{DM}  \end{array} \end{align*}
\end{definition}

The following extends Rieffel's definition \cite{rif}.

\begin{definition}\label{1.3.d} 
Let $\cl{A}, \cl{B}$ be unital dual operator algebras and 
$\cl{F}: \:_{\cl{A}}\mathfrak{DM}\rightarrow \;_{\cl{B}}\mathfrak{DM}$  
be a functor. We say that $\cl{F}$ is a $*$-\textbf{functor} if for 
every pair of objects  $H_1,H_2$ of $_{\cl{A}}\mathfrak{DM}$ 
every operator $F\in \mathrm{ Hom}_{\cl{A}}^{\mathfrak{D}}(H_1,H_2)$ satisfies 
$\cl{F}(F^*)=\cl{F}(F)^*$. 
\end{definition}

\begin{definition}\label{1.4.d} Let $\cl{A}, \cl{B}$ be unital dual operator algebras. 
If there exists an equivalence functor $\cl{F}: \;_{\cl{A}}\mathfrak{M}  
\rightarrow \;_{\cl{B}}\mathfrak{M}  $ which has a $\Delta$-extension as a 
$*$-functor implementing an equivalence between the categories $_{\cl{A}}\mathfrak{DM}  
,\;_{\cl{B}}\mathfrak{DM},$ we say that $\cl{A}$ and $\cl{B}$ are 
\textbf{$\Delta $-equivalent} algebras.
\end{definition}

 In \cite{rif} two $W^*$ algebras $\cl{A}, \cl{B}$ are called Morita equivalent 
if there exists an equivalence of  $\;_{\cl{A}}\mathfrak{M}$ with $\;_{\cl{B}}\mathfrak{M}$ 
implemented by $*$-functors.  
The main theorem of Rieffel for Morita equivalence of $W^*$-algebras
can be formulated as follows \cite[8.5.38]{bm}:

\begin{theorem}\label{1.5} Two $W^*$ algebras $\cl{A}, \cl{B}$ 
are Morita equivalent if and only if they have faithful normal representations 
$\alpha, \beta $ on Hilbert spaces such that the algebras $\alpha(\cl{A}), 
\beta(\cl{B})$ are TRO equivalent.
\end{theorem}

We will generalize this to dual operator algebras:
 
\begin{theorem}\label{1.6} Two unital dual operator algebras $\cl{A}, \cl{B}$ are 
$\Delta $-equivalent if and only if they have completely isometric normal 
representations $\alpha, \beta$ on Hilbert spaces such that the algebras  
$\alpha(\cl{A}), \beta(\cl{B})$ are TRO equivalent.
\end{theorem}

\noindent\emph{For the proof,} we use a recent result obtained jointly 
with V.I. Paulsen \cite{elepaul} (see the Concluding Remarks): If the unital dual operator algebras $\cl{A}, \cl{B}$ have completely isometric normal 
representations with TRO equivalent images then they are stably isomorphic, 
i.e. there exists a Hilbert space $H$ such that the algebras $\cl{A} \overline{\otimes }B(H)$ and $\cl{B}\overline{\otimes }B(H)$ (where $\overline{\otimes }$ denotes the normal spatial tensor product \cite{bm}) 
 are isomorphic as dual operator algebras. One easily checks that the algebras 
$\cl{A}$ and $\cl{A} \overline{\otimes }B(H)$ (resp. $\cl{B}$ and $\cl{B}\overline{\otimes }B(H)$) 
are $\Delta $-equivalent. 

For the converse direction of the proof we need some definitions and facts from \cite{bs}. Let $\cl{A}$ be a unital dual operator algebra. If $K\subset H$ are objects of $_{\cl{A}}\mathfrak{M}$, we say that $K$ is $\cl{A}$-complemented in $H$
if the projection of $H$ onto $K$ belongs to the space $\mathrm{ Hom}_{\cl{A}}(H,H)$. 
We say that the object $H$ is \textbf{$\cl{A}$-universal} if every object $K$ of 
$_{\cl{A}}\mathfrak{M}$ is $_{\cl{A}}\mathfrak{M}$-isomorphic to an  $\cl{A}$-complemented object 
in a direct sum of copies of $H.$ In \cite{bs} it is proved that there exist 
$\cl{A}$-universal objects and that if $(H,\alpha )$ is an $\cl{A}$-universal object then 
$\alpha $ is a complete isometry and $\alpha(\cl{A})= \alpha(\cl{A})^{\prime\prime}.$ Also 
it is proved that there exists a $W^*$ algebra $W^*(\cl{A})$ and a $w^*$-continuous 
completely isometric homomorphism $j: \cl{A}\rightarrow W^*(\cl{A})$ whose range 
generates $W^*(\cl{A})$ as a $W^*$ algebra and which possesses the following universal 
property: given any normal representation $\alpha : \cl{A}\rightarrow B(H),$ there exists 
a unique normal $*$-representation $\overset {\sim }{\alpha }: W^*(\cl{A})\rightarrow B(H)$ 
extending $\alpha.$ An object $H$ is $\cl{A}$-universal if and only if it is 
$W^*(\cl{A})$-universal. 

We now fix unital dual operator algebras $\cl{A}, \cl{B}$ and an equivalence functor $\cl{F}: \;_{\cl{A}}\mathfrak{M}\rightarrow \;_{\cl{B}}\mathfrak{M}$ which has a $\Delta$-extension as a $*$-functor implementing an equivalence between the categories $_{\cl{A}}\mathfrak{DM}$ and $_{\cl{B}}\mathfrak{DM}.$ We still denote the $\Delta$-extension of this functor by $\cl{F}.$ We need the following lemma.

\begin{lemma}\label{restrict} The functor $\cl{F}$ restricts to an equivalence $*$-functor between the categories  $_{W^*(\cl{A})}\mathfrak{M}$ and $_{W^*(\cl{B})}\mathfrak{M}.$
\end{lemma}
\textbf{Proof} If $T\in \mathrm{ Hom}_{W^*(\cl{A})}(H_1,H_2),$ using the fact that $W^*(\cl{A})$ 
(resp. $W^*(\cl{B})$) is a $W^*$-algebra generated by a copy of $\cl{A}$ (resp. $\cl{B}$) and $\cl{F}$ 
is a $*$-functor we can check that $\cl{F}(T)\in \mathrm{ Hom}_{W^*(\cl{B})}(\cl{F}(H_1),\cl{F}(H_2)).$ 
Since the objects of $_{W^*(\cl{A})}\mathfrak{M}$ and $_{\cl{A}}\mathfrak{M}$ coincide, as do the objects of $_{W^*(\cl{B})}\mathfrak{M}$ and $_{\cl{B}}\mathfrak{M}$, we can define 
a functor $\cl{G}:\; _{W^*(\cl{A})}\mathfrak{M}\rightarrow\; _{W^*(\cl{B})}\mathfrak{M}$ by sending 
every object $K$ to the object $\cl{F}(K)$ and every homomorphism $T$ to the 
homomorphism $\cl{F}(T).$ Clearly $\cl{G}$ is a $*$-functor. For every $H_1, H_2\in 
\;_{W^*(\cl{A})}\mathfrak{M}$ the map $\cl{G}: \mathrm{ Hom}_{W^*(\cl{A})}(H_1,H_2)\rightarrow 
\mathrm{ Hom}_{W^*(\cl{B})}(\cl{F}(H_1),\cl{F}(H_2))$ is faithful, being a restriction of $\cl{F}.$ 
Also it is onto because for every $S\in \mathrm{ Hom}_{W^*(\cl{B})}(\cl{F}(H_1),\cl{F}(H_2))$ 
 we can check that $\cl{F}^{-1}(S)\in \mathrm{ Hom}_{W^*(\cl{A})}(H_1,H_2).$ If $K\in \;_{W^*(\cl{B})}
\mathfrak{M},$ since $\cl{F}:\; _{\cl{A}}\mathfrak{M}\rightarrow\; 
_{\cl{B}}\mathfrak{M}$ is an equivalence functor, there exists an object $H\in \;_{\cl{A}}\mathfrak{M}
$ and a unitary $U\in \mathrm{ Hom}_{\cl{B}}(\cl{F}(H),K).$ We can easily 
check that $U$ belongs to $\mathrm{ Hom}_{W^*(\cl{B})}(\cl{F}(H),K).$ It follows that $\cl{G}$ 
is an equivalence $*$-functor. See for example \cite[Theorem 1, section IV-4]{sau}.
$\qquad  \Box$

\begin{corollary}\label{1.7} If $H$ is an $\cl{A}$-universal object then $\cl{F}(H)$ is a $\cl{B}$-universal object.
\end{corollary}
\textbf{Proof} Let $\cl{G}:\; _{W^*(\cl{A})}\mathfrak{M}\rightarrow\; 
_{W^*(\cl{B})}\mathfrak{M}$ be the restriction of $\cl{F}$ as in 
Lemma \ref{restrict}. Every $\cl{A}$-universal object $H$ is  
$W^*(\cl{A})$-universal. Since $\cl{G}$ is an equivalence, $\cl{F}(H)$ is a 
$W^*(\cl{B})$-universal object \cite{rif}, hence $\cl{B}$-universal. $\qquad  \Box$

\medskip

\noindent\emph{We now return to the proof of Theorem \ref{1.6}.} Choose an $\cl{A}$-universal object
$(H,\alpha )$ and denote by $(\cl{F}(H),\beta )$ the corresponding object.
By the previous corollary this object is $\cl{B}$-universal. As we remarked in the discussion before Lemma \ref{restrict} the normal representations $\alpha, \beta$ are complete isometries and the algebras $\alpha(\cl{A}), \beta(\cl{B})$ have the double commutant property: $\alpha(\cl{A})= \alpha(\cl{A})^{\prime\prime}, \beta (\cl{B})= \beta (\cl{B})^{\prime\prime}.$
We denote by $\sigma $ the map
\vspace*{-1ex}
$$\cl{F}: \mathrm{ Hom}_{\cl{A}}^{\mathfrak{D}}(H,H)=\alpha (\Delta (\cl{A}))^\prime\rightarrow \beta (\Delta (\cl{B}))^\prime=\mathrm{ Hom}_{\cl{B}}^{\mathfrak{D}}(\cl{F}(H),\cl{F}(H)),$$
that is $\sigma (T)= \cl F(T),\;T\in  \alpha (\Delta (\cl{A}))'.$ Since 
$\cl{F}: \;_{\cl{A}}{\mathfrak{DM}}\rightarrow \;_{\cl{B}}\mathfrak{DM}$ is an equivalence 
$*$-functor this map is a  $*$-isomorphism. By the $\Delta$-extension property $\sigma $ maps the space $\mathrm{ Hom}_{\cl{A}}(H,H)=\alpha (\cl{A})^\prime$ into $\mathrm{ Hom}_{\cl{B}}(\cl{F}(H),\cl{F}(H))=\beta  (\cl{B})^\prime.$ Since $\cl{F}: \;_{\cl{A}}\mathfrak{M}\rightarrow \;_{\cl{B}}\mathfrak{M}$ is an equivalence functor we have $\sigma (\alpha (\cl{A})^\prime)=\beta (\cl{B})^\prime.$ We define the space $$\cl{M}=\{M: MA=\sigma (A)M\;\text{for\;all\;}A\in 
\alpha (\Delta (\cl{A}))^\prime\}.$$
 By Lemma \ref{1.1} this space is an essential TRO. Choose $M,N \in \cl{M}, B\in \cl{B}.$ 
For all $A\in \alpha (\cl{A})^\prime$ we have $M^*\beta (B)NA=M^*\beta (B)\sigma (A)N.$ 
Since $\sigma (A)\in \beta (\cl{B})^\prime$ the last operator equals 
$M^*\sigma (A)\beta (B)N=AM^*\beta (B)N.$ We proved that 
$\cl{M}^*\beta (\cl{B})\cl{M}\subset \alpha (\cl{A}).$ Symmetrically we can prove 
$\cl{M}\alpha (\cl{A})\cl{M}^*\subset \beta  (\cl{B}).$ It follows 
from \cite[2.1]{ele} that $\alpha (\cl{A})\stackrel{\cl{M}}{\sim}\beta (\cl{B}).$

\section{The generated functor}

In this section we fix unital dual operator algebras $\cl{A, B}$ acting on Hilbert spaces 
$H_0, K_0$ respectively which are TRO  equivalent. We are going to construct a 
functor $\cl{F}_{\cl{U}}$ generated by a $\cl{B}, \cl{A}$ bimodule $\cl{U}.$ In section 3 
we shall prove that every functor implementing the equivalence of Theorem \ref{1.6} is 
unitarily equivalent to such a functor  $ \cl{F}_{\cl{U}} $. 

In \cite[2.8]{ele} it is shown that the TRO 
 $\cl{M}\subset B(H_0,K_0)$ implementing the equivalence can be chosen so that 
$[\cl{M}^*\cl{M}]^{-w^*}=\Delta (\cl{A}),\;\;  [\cl{M}\cl{M}^*]^{-w^*}=\Delta (\cl{B}).$
Define $\cl{U}=[\cl{B}\cl{M}]^{-w^*}, \cl{V}=[\cl{M}^*\cl{B}]^{-w^*}.$
One can now check that
 $\cl{U}=[\cl{M}\cl{A}]^{-w^*}, \cl{V}=[\cl{A}\cl{M}^*]^{-w^*}$
and\vspace*{-1ex}
$$\cl{B}\cl{U}\cl{A}\subset \cl{U},\;\; \cl{A}\cl{V}\cl{B}\subset \cl{V},\;\;
[\cl{V}\cl{U}]^{-w^*}=\cl{A},\;\;[\cl{U}\cl{V}]^{-w^*}=\cl{B}.$$
If $n\in \bb N$ and $S\in Ball(M_{n,1}(\cl{V}))$ we define on the algebraic tensor product
$\cl{U}\otimes H$ a sesquilinear form by the formula:
$$\left\langle T_1\otimes x_1, T_2\otimes x_2\right\rangle_S
=\left\langle\alpha(ST_1)x_1, \alpha(ST_2)x_2 \right\rangle_{H^n}.
$$
We write $ \nor{\cdot}_S $ for the associated seminorm and $\cl{L}_S$ 
for its kernel. The completion of 
$((\cl{U}\otimes H)/\cl{L}_S,\nor{\cdot}_S)$ will be denoted by $H_S$ and the symbol 
$\nor{\cdot}_S$ will be used for the norm of $H_S$ as well. 
Let $\pi_S: \cl U\otimes H\to H_S$ be the quotient map. Again on the algebraic tensor product  $\cl{U}\otimes H$ we define the following seminorm:
\begin{align*}
\nor{\sum_{j=1}^m T_j\otimes x_j}_{\cl F_{\cl U}(H)}=
\sup_{S\in Ball(M_{n,1}(\cl{V})),n\in \mathbb{N}}\nor{\sum_{j=1}^m T_j\otimes x_j}_S
\end{align*}

Since the seminorm $\| \cdot\|_S$ satisfies the parallelogram identity for all $S$,
the previous seminorm satisfies the parallelogram identity too.
If $\cl{L}=\{z\in \cl{U}\otimes H: \|z\|_{\cl F_{\cl U}(H)} =0\}$ the space $(\cl{U}\otimes H)/\cl{L}$
is a pre-Hilbert space. We denote its completion by $\cl{F}_{\cl{U}}(H)$
and we use the same symbol $\nor{\cdot}_{\cl F_{\cl U}(H)}$ for the corresponding norm.
We write $\pi:\cl U\otimes H\to \cl F_{\cl U}(H)$ for the quotient map.
The following lemma is essentially due to Paschke; see for example \cite[8.5.23]{bm}.

\begin{lemma}\label{2.1}There exist partial isometries 
$\{W_k, k\in J\}\subset \cl{M}\;\; (\{V_k, k\in I\} \subset \cl{M})$ such that 
$W_k^*W_k\perp W_m^*W_m \;(V_kV_k^*\perp V_mV_m^*)$ for $k\neq m$ and 
$I_{H_0}=\sum_k\oplus W_k^*W_k \;(I_{K_0}=\sum_k\oplus V_kV_k^*).$
\end{lemma}

\medskip

The following proposition says that we can calculate the norm $\|\cdot\|_{\cl{F}_{\cl{U}}(H)}$ 
using only the  operators $\{S: S\in Ball(M_{n,1}(\cl{M^*})),n\in \mathbb{N}\}.$

\begin{proposition}\label{2.2}
If $\sum_{j=1}^m T_j\otimes x_j \in \cl{U}\otimes H$ then
\begin{align*}
\nor{\sum_{j=1}^m T_j\otimes x_j}_{\cl{F}_{\cl{U}}(H)}=\sup_{S\in Ball(M_{n,1}(\cl{M^*})),n\in \mathbb{N}}\nor{\sum_{j=1}^m T_j\otimes x_j}_S
\end{align*}
\end{proposition}
\textbf{Proof} For $\epsilon >0$ there exist $n\in \mathbb{N}$ and $S\in 
Ball(M_{n,1}(\cl{V}))$ such that $$\nor{\sum_{j=1}^m T_j\otimes x_j}_{\cl{F}_
{\cl{U}}(H)}-\epsilon < \nor{\sum_{j=1}^m\alpha (ST_j)x_j}_{H^n}-\frac{\epsilon }{2}.$$
Using Lemma \ref{2.1} and the fact that $\alpha $ is $w^*$-continuous we can find
 partial isometries $\{V_1,...,V_N \}\subset \cl{M}$ such that
the operator $\sum_{i=1}^NV_iV_i^*$ is a projection and
\begin{align*}
&\nor{\sum_{j=1}^m\alpha (ST_j)x_j}_{H^n}-\frac{\epsilon }{2} \leq \nor{\sum_{j=1}^m
\alpha \left(S\sum_{k=1}^NV_kV_k^* T_j\right)x_j}_{H^n}\\
=&\nor{\alpha \left(S(V_1,...,V_N)\right)\sum_{j=1}^m \alpha ((V_1^*,...,V_N^*)^t T_j) x_j}_{H^n}
\end{align*}
Observe that
$(V_1^*,...,V_N^*)^t$ is in $ Ball(M_{N,1}(\cl{M^*}))$.
So since $\alpha $ is a complete contraction we have
\begin{align*}
\nor{\sum_{j=1}^m T_j\otimes x_j}_{\cl{F}_{\cl{U}}(H)}-\epsilon 
\leq \sup_{S\in Ball(M_{r,1}(\cl{M^*})),r\in \mathbb{N}}\nor{\sum_{j=1}^m T_j\otimes x_j}_S
\end{align*}
Since $\epsilon $ is arbitrary the proof is complete. $\qquad \Box$

\bigskip

For all $n\in \mathbb{N}$ and $S\in Ball(M_{n,1}(\cl{V}))$ we have
$\nor{\pi_s(\xi)}_S \le \nor{\pi(\xi)}_{\cl F_{\cl U}(H)}$ for every $\xi\in \cl U\otimes H.$
This shows that
the map $\pi(\xi)\to \pi_S(\xi)$
is well defined and extends to a contraction $\theta _S: \cl{F}_{\cl{U}}(H)\rightarrow H_S$
between the associated completions.

\begin{lemma}\label{2.3} If $\theta _S: \cl{F}_{\cl{U}}(H)\rightarrow H_S$ is given by 
$\theta_S(\pi(\xi))=\pi_S(\xi)$, 
$$\cl{F}_{\cl{U}}(H)=[\theta _S^*(\pi_S(T\otimes x)): 
S\in Ball(M_{m,1}(\cl{V})), m\in \mathbb{N}, T\in \cl{U}, x\in H]^-$$
\end{lemma}
\textbf{Proof} Let $z\in \cl{F}_{\cl{U}}(H)$ be such that $\left\langle \theta _S^*(\pi_S(T\otimes x)), z\right\rangle_{\cl{F}_{\cl{U}}(H)}=0$ for all $m\in \mathbb{N}, S\in Ball(M_{m,1}(\cl{V})), T\in \cl{U}$ and $x\in H.$ Then $\left\langle \pi _S(T\otimes x), \theta _S(z)\right\rangle_{H_S}=0$ for all $m\in \mathbb{N}, S\in Ball(M_{m,1}(\cl{V})), T\in \cl{U}$ and $x\in H.$ It follows that $\theta _S(z)=0$ for all $m\in \mathbb{N}, S\in Ball(M_{m,1}(\cl{V}))$. But
\begin{equation*}
\|z\|_{\cl{F}_{\cl{U}}(H)}=\sup_{S\in Ball(M_{n,1}(\cl{V})),n\in \mathbb{N}}\|\theta _S(z)\|_S.
\end{equation*}
Indeed, this holds when $z\in \pi(\cl U\otimes H)$ and
it is a standard fact that the equality extends to all $z\in \cl{F}_{\cl{U}}(H)$.
It follows that $z=0.\qquad  \Box$

\medskip

We will show below that the space  $\pi(\cl M\otimes H)$ is dense in $\cl F_{\cl U}(H)$. 
In fact 
we shall prove the following stronger result:

\begin{lemma}\label{2.4}
Let $L$ be an invariant projection for $\alpha (\cl{A}).$ If $T\in \cl{U}$ and $x\in H$  then
$$\pi(T\otimes L(x))\in [\pi(N\otimes L(y)): N\in \cl{M}, y\in H]^{-\cl{F}_{\cl{U}}(H)}.$$
\end{lemma}
\textbf{Proof} On the algebraic tensor product  $\cl{M}^*\otimes
\cl{U}\otimes L(H)$ we define the following sesquilinear form 
$$\left\langle M_1^*\otimes T_1\otimes L(x_1),M_2^*\otimes T_2\otimes L(x_2)\right\rangle
=\left\langle \alpha (M_1^*T_1)L(x_1), \alpha (M_2^*T_2)L(x_2) \right\rangle_H.$$
If $\cl{K}$ is the kernel of $\left\langle \cdot, \cdot \right\rangle$ we denote by 
$K$ the completion of $(\cl{M}^*\otimes \cl{U}\otimes L(H))/\cl{K}$ under the 
corresponding norm and by $\pi_K$ the quotient map
$\cl{M}^*\otimes \cl{U}\otimes L(H)\to K.$ Since the identity operator belongs to 
$[\cl{M}^*\cl{M}]^{-w^*}$ and $\alpha $ is $w^*$-continuous we can check that the space 
$K_1$ generated by vectors of the form $\pi_K(M^*\otimes N\otimes L(y))$ where  
$ M,N \in \cl{M}, y\in H$ is dense in $K.$

\medskip

\noindent\emph{Claim}
For every $N, M_0\in \cl{M}, T\in \cl{U}$ and $x\in H,$ 
$$\pi(NM_0^*T\otimes L(x))\in [\pi(M\otimes L(y)): M\in \cl{M}, y\in H]^{-\cl{F}_{\cl{U}}(H)}.$$

\noindent\emph{Proof} For every
$n\in \bb{N}, S\in Ball(M_{n,1}(\cl{V})), M_i\in \cl{M}, T_i\in \cl{U}, x_i\in H, i=1,...,m$
and $N\in Ball(\cl{M})$ we have
\begin{align*} &\nor{\alpha \left( SN\sum_{i=1}^mM_i^*T_i\right)L(x_i)}_{H^n}
=\nor{\alpha (SN)\sum_{i=1}^m\alpha (M_i^*T_i)L(x_i)}_{H^n}\\
&\leq \nor{\sum_{i=1}^m\alpha (M_i^*T_i)L(x_i)}_{H}=
\nor{\pi_K\left(\sum_{i=1}^m M_i^*\otimes T_i\otimes L(x_i)\right)}_K.
\end{align*}
 It follows from the definition of $\nor{\cdot}_{\cl{F}_{\cl{U}}(H)}$ that
\begin{equation}\label{in21}
\nor{\sum_{i=1}^m NM_i^*T_i\otimes L(x_i)}_{\cl{F}_{\cl{U}}(H)}\leq
\nor{\pi_K\left(\sum_{i=1}^m M_i^*\otimes T_i\otimes L(x_i)\right)}_K.
\end{equation}

Now fix $N\in Ball(\cl{M}), M_0\in \cl{M}, T\in \cl{U}, x\in H$ and $\epsilon >0.$
By the density of $K_1$ in $K$ there exist $N_i,M_i\in \cl{M}, x_i\in H, i=1,...m$ such that
$$\nor{\pi_K(M_0^*\otimes T \otimes L(x))-\pi_K\left(\sum_{i=1}^m M_i^*\otimes N_i\otimes 
L(x_i)\right)}_K < \epsilon.$$
It follows from (\ref{in21}) that
$$\nor{NM_0^*T \otimes L(x)-\sum_{i=1}^m NM_i^*N_i\otimes L(x_i)}_{\cl{F}_{\cl{U}}(H)}
< \epsilon.$$
This proves the Claim. Let $T\in \cl{U}$ and $x\in H.$ It now suffices to show that
$$\label{span} \pi(T\otimes L(x))\in [\pi(NM^*U\otimes L(y)): N, M \in 
\cl{M}, U\in \cl{U}, y\in H]^{-\cl{F}_{\cl{U}}(H)}.$$
Recall the partial isometries $\{V_k, k\in I\}\subset \cl{M}$ from Lemma \ref{2.1}.
We have
\begin{align*}
\lim_{E\subset I,\text{finite}}
&\sca{ \pi(T\otimes L(x))-\sum_{k\in E}\pi(V_kV_k^*T\otimes L(x)) ,
\theta _S^*(\pi_S(U\otimes y))}_{\cl{F}_{\cl{U}}(H)}\\
=\lim_E&\sca{\theta _S\left(\pi\left(T\otimes L(x)-\sum_{k\in E}V_kV_k^*T\otimes 
L(x)\right)\right),
\pi_S(U\otimes y)}_S\\
=\lim_E&\sca{\pi_S\left(T\otimes L(x)-\sum_{k\in E}V_kV_k^*T\otimes L(x)\right), 
\pi_S(U\otimes y)}_S\\
=\lim_E&\sca{\alpha (ST) L(x)-\sum_{k\in E}\alpha (SV_kV_k^*T)L(x), \alpha (SU)(y)}_{H^n}\\
=\lim_E&\sca{\alpha (S(I-\sum_{k\in E} V_kV_k^*)T)L(x), \alpha (SU)(y)}_{H^n}=0.
\end{align*}
Since this net is uniformly bounded from Lemma \ref{2.3} the equality $\pi(T\otimes L(x))=\sum_{k\in I}\pi
(V_kV_k^*T\otimes L(x))$ follows and the proof is complete.
$\qquad  \Box$

\begin{corollary}\label{2.5} The subspace $\pi (\cl{M}\otimes H)$ of $\cl{F}_{\cl{U}}(H)$ is dense.
\end{corollary}

We define a map $\beta : \cl{B}\rightarrow B(\cl{F}_{\cl{U}}(H))$
given by
$$ \beta
(B)(\pi(T\otimes  x))=  \pi(BT\otimes x),\quad  B\in \cl{B},\, T\in \cl{U},\, x\in H.$$
This is a well-defined unital algebraic homomorphism and a contraction. We shall prove the following stronger
result.

\begin{proposition}\label{3.1} The map $\beta $ is a complete contraction.
\end{proposition}
\textbf{Proof} Let $n\in \mathbb{N}$ and $(B_{ij})\in M_n(\cl{B}).$ Fix 
vectors $z_j=\sum_{i=1}^{k_j}\pi (T_i^j\otimes x_i^j),$\\ $ j=1,...,n$ of 
the space $\cl{F}_{\cl{U}}(H)$ and denote by $z$ the vector 
$(z_1,...,z_n)^t.$ Also write $y=\beta ((B_{ij}))(z).$ Then $$\|y\|^2=\sum_{k=1}^n \nor{ \sum_{j=1}^n\sum_{i=1}^{k_j}B_{kj}T_i^j\otimes x_i^j }^2
_{\cl{F}_{\cl{U}}(H)}.$$
 By the definition of the norm of the space $\cl{F}_{\cl{U}}(H),$ given 
$\epsilon >0$ there exist $r\in \mathbb{N},\; S_k\in Ball(M_{r,1}(\cl{V})),\;k=1,...,n$ 
such that
$$ \|y\|^2-\epsilon \leq  \sum_{k=1}^n \nor
{ \sum_{j=1}^n\sum_{i=1}^{k_j} \alpha (S_kB_{kj} T_i^j)(x^j_i) }^2_{H^r}-\frac{\epsilon }{2}
.$$
Since $\alpha $ is $w^*$-continuous from Lemma \ref{2.1} we can find partial isometries 
$V_1,...,V_N\in \cl{M}$ such that
$$ \sum_{k=1}^n \nor
{ \sum_{j=1}^n\sum_{i=1}^{k_j} \alpha (S_kB_{kj} T_i^j)(x^j_i) }^2_{H^r}-\frac{\epsilon }{2}
\leq \sum_{k=1}^n\nor{\sum_{j=1}^n\sum_{i=1}^{k_j} \alpha \left(S_kB_{kj}\sum_{l=1}^N
V_lV_l^*T_i^j\right)(x^j_i)}^2.$$
Let $V=(V_1,...,V_N).$ Now $\alpha $ is an algebraic homomorphism, and hence 

\begin{align*} \|y\|^2&-\epsilon \leq 
\nor{ \alpha ((S_iB_{ij}V)_{1\leq i,j\leq n})\left[\begin{array}{clr} \sum_{i=1}^
{k_1} \alpha (V^*T_i^1)(x_i^1) \\ \vdots \\ \sum_{i=1}^{k_n} \alpha (V^*T_i^n)(x_i^n)   
\end{array}\right]}^2 \end{align*}

Since $(S_iB_{ij}V)_{1\leq i,j\leq n}=(S_1\oplus ...\oplus S_n)(B_{ij})
(V\oplus ...\oplus V)$ and $\|(S_1\oplus ...\oplus S_n)\|\leq 1,$ 
$\|(V\oplus ...\oplus V)\| \leq 1 $ it follows that 
$\|(\alpha (S_iB_{ij}V))\|\leq \|(B_{ij})\|$ and hence
\begin{align*}&\|y\|^2-\epsilon \leq \|(B_{ij})\|^2\sum_{j=1}^n\nor{\sum_{i=1}^{k_j} \alpha 
(V^*T_i^j)(x_i^j)}^2 \\& \leq  \sum_{j=1}^n \nor{ \sum_{i=1}^{k_j}T_i^j
\otimes x_i^j}^2_{\cl{F}_{\cl{U}}(H)}\|(B_{ij})\|^2=\|z\|^2\|(B_{ij})\|^2.
\end{align*}
But $\epsilon$ is arbitrary and so $\nor{\beta((B_{ij}))(z)}=\nor{y}\le \nor{(B_{ij})}\nor{z},$ 
hence $\nor{\beta((B_{ij}))}$ $\le \nor{(B_{ij})}$. Since $n$ is arbitrary,
this shows that $\beta$ is a complete contraction.$\Box$

\begin{proposition}\label{3.2} The map $\beta $ is $w^*$-continuous.
\end{proposition}
\textbf{Proof} Since $\beta $ is a bounded map it suffices to show that given a net
$(B_i)\subset Ball(\cl{B})$ which converges to $0$ in the weak operator topology, the net $(\beta (B_i))$
also converges to $0$ in the weak operator topology. Indeed, for all 
$T_1,T_2\in \cl{U}, x_1,x_2\in H$, $n\in \mathbb{N}$ and $S\in Ball(M_{n,1}(\cl{V})):$
\begin{align*}
&\sca{\beta (B_i)(\pi(T_1\otimes x_1)), \theta_S^*(\pi_S(T_2\otimes x_2))} _{\cl{F}_{\cl{U}}(H)}\\
=&\sca{\theta _S(\pi((B_iT_1\otimes x_1)),\pi_S( T_2\otimes x_2)}_{H_S}
=\sca{\alpha (SB_iT_1)(x_1), \alpha (ST_2)(x_2)}\rightarrow 0. 
\end{align*} The conclusion follows 
from Lemma \ref{2.3}.
$\Box$

\medskip

In the rest of this section if $H\in \,_{\cl{A}}\mathfrak{M}$ we identify $\cl{U}\otimes H$ 
with its image in $\cl{F}_{\cl{U}}(H).$ From the above discussion we have a correspondence $H\in \,_{\cl{A}}\mathfrak{M} \rightarrow 
\cl{F}_{\cl{U}}(H)\in \, _{\cl{B}}\mathfrak{M} .$ If $(H_i,\alpha _i)\in \,_{\cl{A}}\mathfrak{M}, 
i=1,2$ we define a map $\cl{F}_{\cl{U}}(F)$ from the space 
$\cl{F}_{\cl{U}}(H_1)$ into the space $\cl{F}_{\cl{U}}(H_2)$ by the formula  
$$\cl{F}_{\cl{U}}(F)(T\otimes x)=T\otimes F(x)\;\text{for\;all\;}T\in \cl{U}, x\in H_1.$$ 
We can easily check that this map is bounded with norm at most $\|F\|$ and  
$\cl{F}_{\cl{U}}(F)\in \mathrm{ Hom}_{\cl{B}}(\cl{F}_{\cl{U}}(H_1), \cl{F}_{\cl{U}}(H_2)).$ 
This definition completes the definition of the functor $\cl{F}_{\cl{U}}: \;_{\cl{A}}
\mathfrak{M}\rightarrow  \;_{\cl{B}}\mathfrak{M}.$ 

\begin{theorem}\label{4.1}The functor $\cl{F}_{\cl{U}}$ has a $\Delta$-extension.
\end{theorem}
\textbf{Proof} Let $F\in \mathrm{ Hom}_{\cl{A}}^\mathfrak{D}
(H_1,H_2).$ Suppose that $M_1,...,M_m\in \cl{M}$ and $x_1,...,x_m$ $\in H.$ 
 If $n\in \mathbb{N}$ and $S\in Ball(M_{n,1}(\cl{M}^*))$ we have:
\begin{align*}&\nor{\sum_{i=1}^m\alpha _2(SM_i)F(x_i)}=\nor{F^{(n)} 
\sum_{i=1}^m\alpha _1(SM_i)(x_i)}\qquad (F^{(n)}=(F\oplus F\oplus ...\oplus F)) \\\leq &\|F\|\nor{\sum_{i=1}^m\alpha _1(SM_i)
(x_i)}\leq \|F\| \nor{\sum_{i=1}^mM_i\otimes x_i}_{\cl{F}_{\cl{U}}(H_1)}.\end{align*} 
From Proposition \ref{2.2} it follows that $$\nor{\sum_{i=1}^m M_i\otimes F(x_i)}_{\cl{F}_{\cl{U}}(H_2)}\leq \|F\| \nor{\sum_{i=1}^m M_i\otimes x_i}_{\cl{F}_{\cl{U}}(H_1)}.$$
 So we can define a map $\delta (F)$ from the subspace $\cl{M}\otimes H_1$ of $\cl{F}_{\cl{U}}(H_1)
$ into the space $\cl{F}_{\cl{U}}(H_2)$ by the formula  
\begin{equation}\label{eq41}\delta (F)(M\otimes x)=M\otimes F(x)\;\text{for\;all\;}
M\in \cl{M}, x\in H_1. \end{equation}
The map $\delta (F)$ is bounded with norm at most $\|F\|.$ By Corollary \ref{2.5} the space 
$\cl{M}\otimes H_1$ is dense in $\cl{F}_{\cl{U}}(H_1), $ so this map extends to 
$\cl{F}_{\cl{U}}(H_1).$ Since $\Delta (\cl{B})\cl{M}\subset \cl{M},$ equality (\ref{eq41}) 
shows 
that $\delta (F)\in \mathrm{ Hom}_{\cl{B}}^{\mathfrak{D}}( \cl{F}_{\cl{U}}(H_1) ,\cl{F}_{\cl{U}}(H_2)).$ 
  Observe that if $F\in \mathrm{ Hom}_{\cl{A}}(H_1,H_2)$ then $\cl{F}_{\cl{U}}(F)=\delta (F),$ because both operators are bounded and coincide in the dense subspace $\cl{M}\otimes H_1$ of $\cl{F}_{\cl{U}}(H_1).$ Therefore we may define a functor $_{\cl{A}}\mathfrak{DM}
\rightarrow \;_{\cl{B}}\mathfrak{DM}$ by sending every object $H$ to $\cl{F}_{\cl{U}}(H)$ and every homomorphism $F$ to $\delta (F).$ Clearly this functor is a $\Delta$-extension of the functor $\cl{F}_{\cl{U}}.  \qquad \Box$

\medskip

\begin{definition}\label{4.2.d} In the sequel the $\Delta$-extension of the functor 
$\cl{F}_{\cl{U}}$ will be denoted again by $\cl{F}_{\cl{U}}$ and every homomorphism 
$\delta (F)$ defined by equation (\ref{eq41}) by $\cl{F}_{\cl{U}}(F).$
\end{definition}

 Now we will prove that the $\Delta $-extension of $\cl{F}_{\cl{U}}$ is a $*$-functor.

\begin{lemma}\label{4.2} If $U\in \mathrm{ Hom}_{\cl{A}}^\mathfrak{D}
(H_1,H_2)$ is a partial isometry then $$\cl{F}_{\cl{U}}(U^*)=\cl{F}_{\cl{U}}(U)^*.$$
\end{lemma}
\textbf{Proof} Let $M_j\in \cl{M}, x_j\in H_1, 1\leq j \leq m, \;S=(N_1^*,...,N_n^*)
^t \in Ball(M_{n,1}(\cl{M}^*)) .$ We have \begin{align*}
&\nor{\sum_{j=1}^m\alpha _1(SM_j)U^*U(x_j)}^2=
\sum_{i=1}^n\nor{\sum_{j=1}^m\alpha _1(N_i^*M_j)U^*U(x_j)}^2\\
=&\sum_{i=1}^n\nor
{U^*\left(U\sum_{j=1}^m\alpha _1(N_i^*M_j)(x_j)\right)}^2=\sum_{i=1}^n\nor
{\sum_{j=1}^mU\alpha _1(N_i^*M_j)(x_j)}^2 \\=&\sum_{i=1}^n\nor{\sum_{j=1}^m\alpha 
_2(N_i^*M_j)U(x_j)}^2 =\nor{\sum_{j=1}^m\alpha _2
(SM_j)U(x_j)}^2.\end{align*}
Since $S$ was arbitrary in $Ball(M_{n,1}(\cl{M}^*))$ we have $$\nor{\sum_{j=1}^mM_j\otimes 
U^*U(x_j)}_{\cl{F}_{\cl{U}}(H_1)}=\nor{\sum_{j=1}^mM_j\otimes U(x_j)}_
{\cl{F}_{\cl{U}}(H_2)}$$ or equivalently  $$\nor{\cl{F}_{\cl{U}}(U^*U)\left
(\sum_{j=1}^mM_j\otimes x_j\right)}_{\cl{F}_{\cl{U}}(H_1)}=\nor{\cl{F}_{\cl{U}}
(U)\left(\sum_{j=1}^mM_j\otimes x_j\right)}_{\cl{F}_{\cl{U}}(H_2)}.$$
 By Corollary \ref{2.5} we have that \begin{equation}\label{equa}\| \cl{F}_{ \cl{U}}(U^*U)(z)
 \|_{ \cl{F}_{ \cl{U}}(H_1)}= 
\| \cl{F}_{ \cl{U}}(U)(z) \|_{ \cl{F}_{ \cl{U}}(H_2)}\;\;\text{\; for\; all\;\;}z\in 
\cl{F}_{\cl{U}}(H_1).\end{equation}
We proved in Theorem \ref{4.1} that the map $\cl{F}_{\cl{U}}$ between the spaces 
of homomorphisms
is a contraction therefore 
$\cl{F}_{ \cl{U}}(U^*U)$ is an orthogonal projection. It follows now by (\ref{equa}) that
\begin{align*}
\left \langle \cl{F}_{ \cl{U}}(U^*U)(z), z \right \rangle_{ \cl{F}_
{ \cl{U}}(H_1)}=
\left \langle \cl{F}_{ \cl{U}}(U)^*\cl{F}_{ \cl{U}}(U)(z), z \right 
\rangle_{\cl{F}_{ \cl{U}}(H_1)}
\end{align*}
for all $z\in \cl{F}_{\cl{U}}(H_1)$ and so
$\cl{F}_{\cl{U}}(U^*)\cl{F}_{\cl{U}}(U)=\cl{F}_{\cl{U}}(U^*U)=\cl{F}_
{\cl{U}}(U)^*\cl{F}_{\cl{U}}(U).$ Let $W=\cl{F}_{\cl{U}}(U), V=\cl{F}_{\cl{U}}(U^*).$ 
We have proved that 
$\label{eq42}VW=W^*W.$
 Similarly working with the partial isometry $U^*$ we obtain 
$WV=V^*V.$
Now we have $V=\cl{F}_{\cl{U}}(U^*)=\cl{F}_{\cl{U}}(U^*UU^*)
=VWV.$ It follows that $V=W^*WV\Rightarrow V^*=V^*W^*W=V^*VW=WVW=
\cl{F}_{\cl{U}}(UU^*U)=\cl{F}_{\cl{U}}(U)=W$ or equivalently 
$\cl{F}_{\cl{U}}(U^*)=\cl{F}_{\cl{U}}(U)^*. \qquad \Box$

\begin{theorem}\label{4.3} The functor $\cl{F}_{\cl{U}}: \;_{\cl{A}}\mathfrak{DM}
\rightarrow \;_{\cl{B}}\mathfrak{DM}
$ is a $*$-functor.
\end{theorem}
\textbf{Proof} Let $T\in \mathrm{ Hom}_{\cl{A}}^\mathfrak{D}(H_1,H_2)$ with polar decomposition
$T=U|T|.$ Observe that 
$(|T|+\epsilon I)^{-1} \in \alpha _1(\Delta (\cl{A}))^\prime$ for every $\epsilon >0.$
Since $U=w^*-\underset{\epsilon \rightarrow 0}{\lim}\; T(|T|+\epsilon I)^{-1}$ it follows that
$U\in \mathrm{ Hom}_{\cl{A}}^\mathfrak{D}  (H_1,H_2).$ The map
$$\cl{F}_{\cl{U}}: \mathrm{ Hom}_{\cl{A}}^\mathfrak{D}(H_1,H_1)=\alpha _1(\Delta (\cl{A}))'
\to \beta_1(\Delta (\cl{B}))'=\mathrm{ Hom}_{\cl{B}}^\mathfrak{D}(\cl{F}_{\cl{U}}(H_1), \cl{F}_{\cl{U}}(H_1))$$
is an algebraic homomorphism between von Neumann algebras. We also proved in Theorem \ref{4.1}
that it is a contraction. It follows that it is a $*$-homomorphism.
Therefore $\cl{F}_{\cl{U}}(|T|)\geq 0.$ Using the previous lemma we obtain
\begin{align*} \cl{F}_{\cl{U}}(T^*)=&\cl{F}_{\cl{U}}(|T|U^*)=
\cl{F}_{\cl{U}}(|T|)\cl{F}_{\cl{U}}(U^*)\\
=&\cl{F}_{\cl{U}}(|T|)\cl{F}_{\cl{U}}(U)^*
=(\cl{F}_{\cl{U}}(U)\cl{F}_{\cl{U}}(|T|))^*=\cl{F}_{\cl{U}}(T)^* \qquad \Box
\end{align*}

\section{Equivalence functors}

In this section we prove that every functor $\cl{F}$ implementing the
equivalence of Theorem \ref{1.6} is equivalent to a functor of the form $\cl{F}_{\cl{U}}$
for some $\cl B$, $\cl A$ bimodule $\cl{U}$ and we also prove that $\cl{F}$ is normal and completely isometric.

 Throughout this section we fix unital dual operator algebras $\cl{A}, \cl{B}$ and a functor $\cl{F}$
 implementing the equivalence of Theorem \ref{1.6}. We choose an $\cl{A}$-universal object 
$(H_0,\alpha _0)$ . Suppose that $(\cl{F}(H_0),\beta _0)$ is the corresponding object which is 
$\cl{B}$-universal (Corollary \ref{1.7}.) By the proof of Theorem \ref{1.6} (section 1) the map $$\cl{F}: \mathrm{ Hom}_{\cl{A}}^\mathfrak{D}
(H_0,H_0)=\alpha _0(\Delta (\cl{A}))^\prime\rightarrow \beta _0(\Delta (\cl{B}))^\prime=\mathrm{ Hom}_{\cl{B}}^\mathfrak{D}
(\cl{F}(H_0),\cl{F}(H_0))$$ is a $*$-isomorphism with the property $\cl{F}(\alpha _0(\cl{A})^\prime)=\beta _0(\cl{B})^\prime,$ the space
$$\cl{M}=\{M\in B(H_0,\cl{F}(H_0)): MF=\cl{F}(F)M\;\text{for\;all\;}F\in \alpha _0(\Delta (\cl{A}))^\prime\}$$
is an essential TRO and the algebras $\alpha_0(\cl{A}), \beta_0(\cl{B})$ are TRO equivalent
via the space $\cl{M}.$ We denote by $\cl U$ and $\cl V$ the spaces
$$\cl{U}=[\cl{M}\alpha _0(\cl{A})]^{-w^*}, \quad \cl{V}=[\alpha _0(\cl{A})\cl{M}^*]^{-w^*}$$
which satisfy the following relations 
$$\beta _0(\cl{B})\cl{U}\alpha _0(\cl{A})\subset \cl{U},
\;\; \alpha _0(\cl{A})\cl{V}\beta _0(\cl{B})\subset \cl{V},
\;\;[\cl{V}\cl{U}]^{-w^*}=\alpha _0(\cl{A}),
\;\;[\cl{U}\cl{V}]^{-w^*}=\beta _0(\cl{B}).$$
As in section 2 we define a functor $\cl{F}_{\cl{U}}: \;_{\cl{A}}\mathfrak{M}\rightarrow 
 \;_{\cl{B}}\mathfrak{M} $ 
which has a $\Delta$-extension.
In the rest of this section for every $(H,\alpha )\in \;_{\cl{A}}\mathfrak{M}$
 we identify the element $T\otimes x$ with its image in  $\cl{F}_{\cl{U}}(H)$ (see section 2). 
Also we identify the algebra $\alpha _0(\cl{A})$ with $\cl{A}$ and 
the algebra $\beta _0(\cl{B})$ with $\cl{B}.$ 

\begin{lemma}\label{6.1}
(i) The map $T\otimes x \to T(x)\quad T\in \cl{U},\, x\in H_0$ extends to a unitary
$U: \cl{F}_{\cl{U}}(H_0)\rightarrow \cl{F}(H_0)$
which belongs to the space $\mathrm{ Hom}_{\cl B}(\cl{F}_{\cl U}(H_0), \cl F(H_0)).$

(ii) For all $F\in \mathrm{ Hom}_{\cl{A}}^\mathfrak{D}  (H_0,H_0)$ the equality 
$U\cl{F}_{\cl{U}}(F)=\cl{F}(F)U$ holds.
\end{lemma}
\textbf{Proof} (i) For all $T_1,...,T_m\in \cl{U}, x_1,...,x_m\in H_0$ we have $$\nor{\sum_{j=1}^m T_j\otimes x_j}_{\cl{F}_{\cl{U}}(H_0)}=\sup_{S\in Ball(M_{n,1}(\cl{V})),n\in \mathbb{N}}\nor{\sum_{j=1}^m ST_j(x_j)}_{H_0^{(n)}}\leq \nor{\sum_{j=1}^m T_j(x_j)}_{\cl{F}(H_0)}.$$

For arbitrary $\epsilon >0$ there exist (Lemma \ref{2.1}) partial isometries
$V_1,...,V_n\in \cl{M}$ such that the operator $\sum_{i=1}^nV_iV_i^*$ is a projection and
\begin{align*}
& \nor{\sum_{j=1}^m T_j(x_j)}_{\cl{F}(H_0)}-\epsilon
\leq  \nor{ \sum_{l=1}^n V_lV_l^*\sum_{j=1}^m T_j(x_j) }_{\cl{F}(H_0)}\\
\leq & \nor{ (V_1^*,...,V_n^*)^t \sum_{j=1}^m T_j(x_j) }_{H_0^{(n)}}
\leq \nor{\sum_{j=1}^m T_j\otimes x_j}_{\cl{F}_{\cl{U}}(H_0)}.
\end{align*}
It follows that
$\nor{\sum_{j=1}^m T_j\otimes x_j}_{\cl{F}_{\cl{U}}(H_0)}=\nor{\sum_{j=1}^m T_j(x_j)}_{\cl{F}(H_0)}.$
So the map $T\otimes x\rightarrow T(x), T\in \cl{U}, x\in H_0$ extends to an isometry
$U: \cl{F}_{\cl{U}}(H_0)\rightarrow \cl{F}(H_0).$ Since 
$[\cl{U}(H_0)]^-=\cl{F}(H_0)$ the image of $U$ is dense in $\cl{F}(H_0),$
so $U$ is a unitary. We can easily check that $U\in  \mathrm{ Hom}_{\cl{B}}(\cl{F}_{\cl{U}}(H_0), \cl{F}(H_0)).$

\medskip

(ii) Let $F\in \mathrm{ Hom}_{\cl{A}}^\mathfrak{D}
(H_0,H_0).$ For every $M\in \cl{M}, x\in H_0$ we have $$(U\cl{F}_{\cl{U}}(F))(M\otimes x)=U(M\otimes F(x))=M(F(x))=\cl{F}(F)M(x)=(\cl{F}(F)U)(M\otimes x).$$
 By Corollary \ref{2.5} it follows that $U\cl{F}_{\cl{U}}(F)=\cl{F}(F)U.\qquad  \Box$

\bigskip

 The following lemma is analogous to \cite[Proposition 4.9]{rif}. The proof is similar, using the
$\Delta$-extension of the functors $\cl F, \cl{F}_{\cl{U}}$.

\begin{lemma}\label{6.2} If $\{H_j: j\in I\}$ are objects of $_{\cl{A}}\mathfrak{M}
,$ then there exist unitaries $W\in \mathrm{ Hom}_{\cl{B}}(\oplus _j\cl{F}(H_j), \cl{F}(\oplus _jH_j)),$ and $V\in \mathrm{ Hom}_{\cl{B}}(\oplus _j\cl{F}_{\cl{U}}(H_j), \cl{F}_{\cl{U}}(\oplus _jH_j)).$
\end{lemma}

\begin{theorem}\label{6.3} The functors $\cl{F}, \cl{F}_{\cl{U}}$ are equivalent as functors between the categories $_{\cl{A}}\mathfrak{M},\; _{\cl{B}}\mathfrak{M}
$ and their $\Delta$-extensions are equivalent as $*$-functors between the categories $_{\cl{A}}\mathfrak{DM}
, \;_{\cl{B}}\mathfrak{DM}
.$
\end{theorem}
\textbf{Proof}
Since $H_0$ is an $\cl{A}$-universal object,
it is also $W^*(\cl A)$-universal (section 1). Therefore, by \cite[Proposition 1.1]{rif}
for every $K \in\; _{\cl{A}}\mathfrak{M}$
there exists a set of indices $J_K,$ projections
$$\{Q_i^K: i\in J_K\}\subset \mathrm{ Hom}_{W^*(\cl A)}(H_0,H_0)\subset \mathrm{ Hom}_{\cl{A}}(H_0,H_0)$$
and a unitary
$$W_K\in \mathrm{ Hom}_{W^*(\cl{A})}(K, \oplus _iQ_i^K(H_0))\subset \mathrm{ Hom}_{\cl{A}}(K, \oplus _iQ_i^K(H_0)).$$
Since the $\Delta$-extensions of $\cl{F}, \cl{F}_{\cl{U}}$ are $*$-functors, the operators
$\cl{F}(W_K), \cl{F}_{\cl{U}}(W_K)$ are unitaries. By Lemma \ref{6.2} we can view
$\cl{F}_{\cl{U}}(W_K)$ as an element
$$\cl{F}_{\cl{U}}(W_K)\in \mathrm{ Hom}_{\cl{B}}(\cl{F}_{\cl{U}}(K), \oplus _i\cl{F}_{\cl{U}}(Q_i^K(H_0)))$$
and
$$\cl{F}(W_K)\in \mathrm{ Hom}_{\cl{B}}(\cl{F}(K), \oplus _i\cl{F}(Q_i^K(H_0))).$$
 Lemma \ref{6.1},ii shows that $U\cl{F}_{\cl{U}}(Q_i^K)=\cl{F}(Q_i^K)U.$ Thus the operator
$$U_i^K=U|_{\cl{F}_{\cl{U}}(Q_i^K(H_0))}: \cl{F}_{\cl{U}}(Q_i^K(H_0))\rightarrow \cl{F}(Q_i^K(H_0))$$
is a unitary for all $i\in J_K.$ So we can define the unitary
$$V_K=\cl{F}(W_K^*)(\oplus _iU_i^K) \cl{F}_{\cl{U}}(W_K)\in \mathrm{ Hom}_{\cl{B}}(\cl{F}_{\cl{U}}(K), 
\cl{F}(K))
\subset \mathrm{ Hom}_{\cl{B}}^{\mathfrak{D}}(\cl{F}_{\cl{U}}(K), \cl{F}(K)).$$
As in the proof of \cite[Proposition 5.4]{rif} we can prove that the unitaries 
$\{V_K: K\in \;_{\cl{A}}\mathfrak{M}   \}$
implement both the required equivalences. $\qquad \Box$

\begin{definition}\label{6.1.d} Let $\cl{A}_1, \cl{B}_1$ be unital dual operator algebras. 
A functor $\cl{G}: \;_{\cl{A}_1}\mathfrak{M}
 \rightarrow \;_{\cl{B}_1}\mathfrak{M}$  is called \textbf{completely isometric (resp. normal)}
if for every pair of objects $H_1, H_2$ the map
$\cl{G}: \mathrm{ Hom}_{\cl{A}_1}(H_1,H_2) \to \mathrm{ Hom}_{\cl{B}_1}(\cl{G}(H_1),\cl{G}(H_2))$
is a complete isometry (resp. $w^*$-continuous). Similarly for a functor $\cl{G}: \;_{\cl{A}_1}\mathfrak{DM}
 \rightarrow \;_{\cl{B}_1}\mathfrak{DM}$.
\end{definition}

\begin{lemma}\label{6.8} The functor $\cl{F}_{\cl{U}}: \;_{\cl{A}}\mathfrak{DM}
 \rightarrow \;_{\cl{B}}\mathfrak{DM}  $ is   normal.
\end{lemma}
\textbf{Proof}  Let $H_1,H_2\in \;_{\cl{A}}\mathfrak{M}  .$ We have proved in Theorem
\ref{4.1} that $\|\cl{F}_{\cl{U}}(F)\|\leq \|F\|$ for all $F\in \mathrm{ Hom}_{\cl{A}}^\mathfrak{D}  (H_1,H_2).$
So it suffices to show that if $(F_i)$ is a bounded 
net of the space $\mathrm{ Hom}_{\cl{A}}(H_1,H_2))$ which converges in the weak
operator topology  to $0$ then the net $(\cl{F}_{\cl{U}}(F_i))$ converges
in the weak operator topology  to $0$ too. We recall from section 2 the contractions 
$\theta _S: \cl{F}_{\cl{U}}(H_2)\rightarrow H_{2,S}$
and the quotient maps $\pi, \pi_S $ where $S\in Ball(M_{n,1}(\cl{V})), n\in \mathbb{N}.$ If 
$M \in \cl{M}, x\in H_1, T\in \cl{U}$ and $y\in H_2$ then 
\begin{align*}&\left\langle \cl{F}_{\cl{U}}(F_i)(\pi (M \otimes x)), 
\theta _S^*(\pi _S(T\otimes y))\right\rangle _{\cl{F}_{\cl{U}}(H_2)}
=\left\langle \theta _S(\pi (M\otimes F_i(x))), \pi _S(T\otimes y)\right
\rangle _{H_{2,S}}\\=&\left\langle \pi _S(M\otimes F_i(x)), \pi _S(T\otimes y)
\right\rangle _{H_{2,S}}
=\left\langle\alpha _2(SM)F_i(x), \alpha _2(ST)(y)\right\rangle \rightarrow 0  \end{align*}

We recall from Lemma \ref{2.3} that
$$\cl{F}_{\cl{U}}(H_2)=
[\theta _S^*(\pi _S(T\otimes y)): S\in Ball(M_{m,1}(\cl{V})), m\in \mathbb{N}, T\in \cl{U}, y\in H_2]^-$$
and from Corollary \ref{2.5} that the space $\pi (\cl{M}\otimes H_1)$ is dense in
$\cl{F}_{\cl{U}}(H_1).$ Since the net $(\cl{F}_{\cl{U}}(F_i))$ is bounded it follows that
$\left\langle \cl{F}_{\cl{U}}(F_i)(z), \xi \right \rangle \rightarrow 0$ for all 
$z\in \cl{F}_{\cl{U}}(H_1), \xi \in \cl{F}_{\cl{U}}(H_2).  \qquad \Box $

\begin{lemma}\label{complete} The functor $\cl{F}_{\cl{U}}: \;_{\cl{A}}\mathfrak{DM}
 \rightarrow \;_{\cl{B}}\mathfrak{DM}  $ is  completely isometric.
\end{lemma}
\textbf{Proof} Let  $(H_1,\alpha _1), (H_2,\alpha _2) \in \;_{\cl{A}}
\mathfrak{M}$ and $(F_{ij})\in M_n( \mathrm{ Hom}_{\cl{A}}^\mathfrak{D}  (H_1,H_2)) $
for $n\in \mathbb{N}.$  Fix vectors
$z_j=\sum_{k=1}^{m_j}M_k^j\otimes x_k^j\in \cl{M}\otimes H_1,\;\; j=1,...,n$
and denote by $z$ the vector $(z_1,...,z_n)^t.$ Then
\begin{align*} \nor{( \cl{F}_{\cl{U}}(F_{ij}) )(z)}^2
=\sum_{i=1}^n\nor{\sum_{j=1}^n\sum_{k=1}^{m_j}M_k^j\otimes F_{ij}(x_k^j)}^2_
{\cl{F}_{\cl{U}}(H_2)}.
\end{align*}
 For $\epsilon >0$ by Proposition \ref{2.2} there exist $r\in \mathbb{N}$ and 
$S_i=(S_1^i,...,S_r^i)^t\in Ball(M_{r,1}(\cl{M}^*)), i=1,...,n$ such that
\begin{align*}& \nor{( \cl{F}_{\cl{U}}(F_{ij}) )(z)}^2 -\epsilon \leq \sum_{i=1}^n
\nor{\sum_{j=1}^n\sum_{k=1}^{m_j}\alpha _2(S_iM_k^j)F_{ij}(x_k^j)}_{H_2^r}^2\\ = & 
\sum_{i=1}^n\sum_{l=1}^r\nor{\sum_{j=1}^n\sum_{k=1}^{m_j}\alpha _2(S^i_lM_k^j)F_{ij}
(x_k^j)}_{H_2}^2=
\sum_{i=1}^n\sum_{l=1}^r
\nor{\sum_{j=1}^n\sum_{k=1}^{m_j}F_{ij}\alpha _1(S^i_lM_k^j)(x_k^j)}_{H_2}^2\\= & \sum_{l=1}^r
 \nor{(F_{ij}) \left[\begin{array}{clr} \sum_{k=1}^{m_1} \alpha_1 (S_l^1M_k^1)(x_k^1) \\ \vdots
\\ \sum_{k=1}^{m_n} \alpha_1 (S_l^nM_k^n)(x_k^n)   \end{array}\right]  }^2_{H_2^n} 
\end{align*}
\begin{align*}\leq &
\|(F_{ij})\|^2\sum_{l=1}^r
 \nor{\left[\begin{array}{clr} \sum_{k=1}^{m_1} \alpha_1 (S_l^1M_k^1)(x_k^1) \\ \vdots
\\ \sum_{k=1}^{m_n} \alpha_1 (S_l^nM_k^n)(x_k^n)   \end{array}\right]  }^2_{H_1^n} =
\|(F_{ij})\|^2\sum_{i=1}^n\nor{\sum_{k=1}^{m_i}\alpha _1(S_iM_k^i)(x_k^i)}^2_{H^r_1}\\
\leq & \nor{(F_{ij})}^2 \sum_{i=1}^n\|z_i\|^2_
{\cl{F}_{\cl{U}}(H_1)} = \nor{(F_{ij})}^2 \|z\|^2_{\cl{F}_{\cl{U}}(H_1)^n}.  \end{align*}

Since $\epsilon $ was arbitrary we have  $ \nor{( \cl{F}_{\cl{U}}(F_{ij}) )(z)}\leq
\|(F_{ij})\|\|z\| $ for all $z\in M_{n,1}(\cl{M}\otimes H_1).$ From Corollary \ref{2.5}
 it follows that $\nor{( \cl{F}_{\cl{U}}(F_{ij}) )}\leq \|(F_{ij})\|.$ 
By Theorem \ref{6.3}, $\cl F_\cl{U}$ is an equivalence functor; hence there is a functor 
$\cl{G}$ such that $\cl{G}\circ\cl{F}_{\cl{U}}$ is equivalent to the identity functor. 
As above we see that $\cl{G}$ can be taken of the form $\cl{G}_{\cl{W}}$ for a suitable bimodule $\cl{W}$. Hence the reverse inequality follows.
$\qquad \Box$

\medskip

Combining Lemmas \ref{6.8}, \ref{complete} and Theorem \ref{6.3} we obtain the next theorem:

\begin{theorem}\label{6.9} Every functor implementing the equivalence of Theorem \ref{1.6}
is a normal and completely isometric functor.
\end{theorem}

\paragraph{Concluding Remarks}\label{rems}\mbox{}
 
\medskip

\noindent\textbf{1. } In a companion paper \cite{ele2} 
we show that every functor implementing the equivalence of Theorem 
\ref{1.6} maps completely isometric representations to  completely isometric representations 
and reflexive algebras to reflexive algebras. Also we present examples of 
$\Delta$-equivalent and $\Delta$-inequivalent CSL algebras. 

\medskip

\noindent\textbf{2. }
The original proof (see \texttt{ArXiv: math.OA/0607489 v.3}) 
of one direction of Theorem \ref{1.6} (if the algebras 
have completely isometric normal representations with 
TRO equivalent images then they are $\Delta $-equivalent) 
was by proving that the functor $\cl{F}_{\cl{U}},$ constructed in section 2, 
 is an equivalence functor. 
After this work was submitted, the present author and V.I. Paulsen 
proved in \cite{elepaul} that TRO equivalent algebras are stably isomorphic. 
We thank the referee for suggesting that we use this result to shorten our original proof.

\medskip

\noindent{\em Acknowledgment:} I would like to express appreciation to
 A. Katavolos for his helpful comments and suggestions during
the preparation of this work, which is part of my doctoral thesis.
Thanks are also due to V. Paulsen for useful discussions. 
This research was partly supported by Special Account Research
Grant No. 70/3/7463 of the University of Athens. 
A preliminary version of these results 
was presented by A. Katavolos
in the Operator Algebra Workshop held at
Queen's University Belfast in May 2006.

\noindent G.K. ELEFTHERAKIS\\
 Department of Mathematics, University of
Athens, \\Panepistimioupolis 157 84, Athens, Greece\\
e-mail: \texttt{gelefth@math.uoa.gr}


\vspace{-1ex}

\begin{thebibliography}{99}
\bibitem{bmor}
Blecher, David P.  A Morita theorem for algebras of operators on Hilbert space.
 {\em J. Pure Appl. Algebra}  156  (2001),  no. 2-3, 153--169.

\bibitem{bmag}
Blecher, David P. ;  Magajna, Bojan. Duality and operator algebras: automatic weak$\sp *$ continuity
and applications.
{\em J. Funct. Anal.}  224  (2005),  no. 2, 386--407.

\bibitem{bm}
Blecher, David P. ;  Le Merdy, Christian. {\em Operator algebras and their modules---an 
operator space approach.}
London Mathematical Society Monographs. New Series, 30. Oxford Science Publications.
The Clarendon Press, Oxford University Press, Oxford,  2004. x+387 pp. ISBN: 0-19-852659-8 

\bibitem{bmp}
Blecher, David P. ;  Muhly, Paul S. ;  Paulsen, Vern I.  Categories of operator modules 
(Morita equivalence and projective
 modules).
 {\em Mem. Amer. Math. Soc.}  143  (2000),  no. 681, viii+94 pp.

\bibitem{bs}
Blecher, David P. ;  Solel, Baruch. A double commutant theorem for operator algebras.
 {\em J. Operator Theory}  51  (2004),  no. 2, 435--453.

\bibitem{ele}
G.K. Eleftherakis, TRO equivalent algebras, preprint, ArXiv:math.OA/0607488.

\bibitem{ele2}
G.K. Eleftherakis, Morita type equivalences and reflexive
algebras, preprint.


\bibitem{elepaul}
G.K. Eleftherakis, V.I. Paulsen, Stably isomorphic dual operator algebras, preprint, 
ArXiv:math.OA/
07052921.

\bibitem{rif}
 Rieffel, Marc A.  Morita equivalence for $C\sp{*} $-algebras and $W\sp{*}
 $-algebras.
 {\em J. Pure Appl. Algebra}  5  (1974), 51--96.

\bibitem{rif2}
Rieffel, Marc A.  Morita equivalence for operator algebras.
 Operator algebras and applications, Part I (Kingston, Ont., 1980), 
 pp. 285--298, {\em Proc. Sympos. Pure Math.,} 38, Amer. Math. Soc., Providence, R.I.,  1982. 


\bibitem{sau}
 MacLane, Saunders . {\em Categories for the working mathematician.}
Graduate Texts in Mathematics, Vol. 5.
Springer-Verlag, New York-Berlin,  1971. ix+262 pp.
\end{thebibliography}
\end{document}